\documentclass{article}
\usepackage[dvips]{graphics}
\usepackage{amsfonts}
\usepackage{amssymb}
\usepackage{amsmath}
\usepackage{amsthm}
\usepackage{latexsym}
\newtheorem{theorem}{Theorem}
\newtheorem{lemma}[theorem]{Lemma}
\newtheorem{proposition}[theorem]{Proposition}
\newtheorem{corollary}[theorem]{Corollary}
\theoremstyle{definition}
\newtheorem{definition}[theorem]{Definition}

\def\Ad{\mathrm{Ad }}
\def\Ad_{\underline{\mathrm{Ad }}}
\def\Kb{K_P}
\def\kb{\mathfrak{k}_P}

\def\Mb{\mathrm{End}(V)_P}
\def\G2{\mathcal{G}^2_{con}}
\def\G4{\mathcal{G}^4_{con}}
\def\Poincare{Poincar\'{e} }
\def\C1{\mathcal{C}^1_{con}}
\def\C3{\mathcal{C}^3_{con}}
\def\k{\mathfrak{k}}

\def\Ck{\mathcal{C}^k_{Dir, A_0}}
\def\k{\mathfrak{k}}
\def\Gk{\mathcal{G}^{k+1}_{Dir}}

\def\Hko{\mathcal{H}^{k}_{Dir, 0}}
\def\Lie{\mathrm{Lie}}
\def\Span{\mathrm{Span}}
\def\tr{\mathrm{tr}}

\begin{document}

\title{On the Holonomy of the Coulomb Connection over Manifolds with 
Boundary}
\author{William E. Gryc
\thanks{\emph{Present address:} Department of Mathematics, Morehouse College, 830 Westview Dr, Atlanta, GA 30314}%
}                     
%
%
%
%
\maketitle
\begin{abstract}
Narasimhan and Ramadas showed in \cite{NR} that the restricted holonomy group of the Coulomb connection is dense in the connected component of the identity of the gauge group when one considers the product principal bundle $S^3\times SU(2)\to S^3$.  Instead of a base manifold $S^3$, we consider here a base manifold of dimension $n\ge 2$ with a boundary and use \emph{Dirichlet boundary conditions} on connections as defined by Marini in \cite{Ma}.  A key step in the method of Narasimhan and Ramadas consisted in showing that the linear space spanned by the curvature form at one
specially chosen connection is dense in the holonomy Lie algebra with respect to an appropriate Sobolev norm.
      Our objective is to explore the effect of the presence of a
boundary on this construction of the holonomy Lie algebra.
Fixing appropriate Sobolev norms, it will be shown that the space
spanned, linearly, by the curvature form at any one connection
is never dense in the holonomy Lie algebra.  In contrast, the linear space
spanned by the curvature form and its first commutators at the flat connection
is dense and, in the $C^\infty$ category, is in fact the entire holonomy
Lie algebra. The former, negative, theorem
is proven for a general principle bundle over $M$, while the latter,
positive, theorem is proven
only for a product bundle over the closure of a bounded open subset of $\mathbb{R}^n$.
Our technique for proving absence of density consists in showing
that the linear space spanned by the
curvature form at one point is contained in the kernel of a linear map consisting of a third order
differential operator, followed by a restriction operation at the boundary; this mapping is determined by the mean curvature
of the boundary.

\end{abstract}
%
%
%
%
%
\section{Introduction}
\label{intro}
In this paper we study the space of connections $\mathcal{A}$ of a certain type of principal bundle and the set of gauge transformations 
$\mathcal{G}$ that act on these connections.  In particular, we are interested in the quotient $\mathcal{A}\to\mathcal{A}/\mathcal{G}$.  This quotient is important in classical Yang-Mills theory; the equivalence classes of $\mathcal{A}/\mathcal{G}$ are physically observable, while individual members of $\mathcal{A}$ are not.  This 
distinction has led to complications in Yang-Mills theory, such as the Gribov ambiguity (see \cite{Gribov}, \cite{NR}, \cite{Singer1}, for example).

In \cite{NR}, Narasimhan and Ramadas considered the Coulomb connection on the quotient $\mathcal{A}^k\to\mathcal{A}^k/\mathcal{G}^{k+1}$, where $\mathcal{A}^k$ and $\mathcal{G}^{k+1}$ are certain Sobolev spaces 
of generic connections and gauge transformations, respectively, of the trivial $SU(2)$ principal bundle over $S^3$ (note that the Coulomb connection is a connection over a space of connections $\mathcal{A}^k$ of the bundle $S^3\times SU(2)\to S^3$).  They showed that in this case the image of the curvature form of the Coulomb connection at the Maurer-Cartan connection is dense in the gauge algebra.  Since the image of the curvature form is contained in the holonomy algebra, this fact implies that the restricted holonomy group of the Coulomb connection is dense in the connected component of the identity of $\mathcal{G}^k$.  They described the physical ramifications of this density, and call it a ``maximal ambiguity'' in reference to the Gribov ambiguity.

Here we are interested in the Coulomb connection when the principal bundle $P$ is over a compact manifold $M$ with 
\emph{non-empty boundary}, and the structure group $K$ is a compact subgroup of $SO(m)$ or $SU(m)$.  This bundle need not 
be trivial.  In this with boundary case, we will want our connections to satisfy \emph{Dirichlet boundary conditions} as 
defined by Marini in \cite{Ma}.  In this case, we see that, unlike in the case found in \cite{NR}, the image of the curvature form of the Coulomb connection cannot be dense in the gauge algebra at any one fixed point.  Indeed, we will show
\\\\
\noindent{\bf Theorem.} \emph{Let $M$ be a compact $n$-dimensional manifold with boundary with $n\ge 2$ and let $\nabla_{A}$ be a 
connection of Sobolev class $k$ for an integral $k>n/2+1$ that satisfies Dirichlet 
boundary conditions.  Define a set $\mathcal{L}_A$ as
\begin{equation*}
\mathcal{L}_A = \Span(\mathrm{Im} (\mathcal{R}_A)),
\end{equation*}
where $\mathcal{R}_A$ is the curvature of the Coulomb connection at $\nabla_A$. Let $\Gk$ be the gauge transformations of Sobolev class $k+1$ that satisfy Dirichlet boundary conditions.
There exists a bounded nonzero operator $T_A:\Lie(\Gk)\to L^2_{k-\frac{5}{2}}(\kb|_{\partial M})$ such that 
$\mathcal{L}_A\subseteq \mathrm{ker}(T_A)$, and thus $\mathcal{L}_A$ cannot be dense in $\Lie(\Gk)$.  This linear map is given by
\begin{equation*}
T_A(f)= d_A\Delta_A f + 2(n-1)H_{}\Delta_A f,
\end{equation*}
where $H_{}$ is the mean curvature of the boundary of $M$.}
\\
\indent The above theorem will be restated as Lemma \ref{bct} and Theorem \ref{bigcor} of Section~\ref{2.3}.  While the image of the curvature form $\mathcal{R}_A$ at a fixed connection $\nabla_A$ may not be dense, the holonomy algebra may still be dense as the span of $\mathcal{R}_A$ for a fixed connection $\nabla_A$ is not the entire holonomy algebra (indeed, this is not even an algebra). The denseness of the gauge algebra was the physically relevent result of Narasimhan and Ramadas, and this denseness will still hold in at least one specific case despite the presence of a boundary and Dirichlet boundary conditions.
\\
\noindent{\bf Theorem.} \emph{Consider the trivial principal bundle $\bar{O}\times K\to\bar{O}$, 
where $O\subseteq\mathbb{R}^n$ is a bounded open set with smooth boundary, $n\ge 2$, and $K$ is a compact subset of $SO(m)$ or $SU(m)$.  The linear space spanned by the curvature form and its first commutators at the flat connection generates all $C^\infty$ elements of the gauge algebra. Furthermore, the restricted holonomy group of 
the Coulomb connection is dense in the connected component of the identity of $\Gk$.}
\\
\indent The above theorem will be restated as Lemma \ref{mainthm} and Theorem \ref{lastcor} in Section~\ref{2.5}.

This paper is a condensed version of the author's doctoral thesis \cite{GrycThesis}.  One difference is that in \cite{GrycThesis} only dimension $3$ is considered, while in this paper we consider any dimension $n$ greater than $2$.   Certain arguments that are 
omitted or shortened in this paper can be found in \cite{GrycThesis} for the specific case of $n=3$, and where appropriate a reference to \cite{GrycThesis} will be given if more detail can be found there. However, please note that some of the notation has been changed between the thesis and this paper.
\section{Background and Notation}
\label{BN}
\subsection{The Differential Geometric Setting of Yang-Mills Theory}\label{1.1}
The differential geometric set-up of Yang-Mills theory can be found in \cite{BL}.  Here we 
give a brief review as well as establish notation.

Let $M$ be a 
compact oriented $n$-dimensional Riemannian manifold with 
boundary for $n\ge 2$, and let $P\to M$ be a principal bundle with a semisimple compact connected 
structure group $K$. We also assume that $K$ acts faithfully on a finite dimensional real 
(or complex) inner product space $V$ by isometries, and 
thus we view $K$ as a compact matrix group and a subgroup of $SO(m)$ (or $SU(m)$, respectively).  
The natural action of $K$ on $V:=\mathbb{R}^m$ (or $V:= 
\mathbb{C}^m$) induces a vector bundle $E:=P\times_{K} V$ (for the definition 
of these associated bundles, see Chapter 1.5 in \cite{KN}).  $K$ also acts on 
itself and its Lie algebra $\mathfrak{k}$ via the adjoint representation, and thus we have the
corresponding bundles $\Kb := P\times_{K} K$ and $\kb := 
P\times_{K} \mathfrak{k}$.

Note that $\kb$ is a vector bundle, while $\Kb$ is not.  However, both $\kb$ 
and $\Kb$ are subbundles of the vector bundle $\Mb:=P\times_K 
\mathrm{End}(V)$, where again $K$ acts by the adjoint action.
 
Recall the exponential map 
$\exp :\mathfrak{k}\to K $.  Since $\Ad_ (k)\circ \exp = \exp\circ
\mathrm{Ad }(k)$, for any $k\in K$, 
we have an induced map $\exp :\kb\to\Kb$.

As $\mathrm{End}(V)$ acts on $V$ in an obvious way, fibers of $\Mb$ act 
on fibers of $E$. Viewing $\kb$ 
and $\Kb$ as subbundles of $\Mb$, fibers of $\kb$ and $\Kb$ also act on 
fibers of $E$. Similar reasoning also tells us that given two elements $\phi,
\psi$ in the same fiber of $\kb$, the Lie bracket $[\phi,\psi]$ is well-defined.

A Koszul connection $\nabla_A$ on $E$ induces a Koszul connection 
also called $\nabla_A^{Hom}$ on $\mathrm{Hom}(E,E)$ 
(see \cite{BL}, \cite{Dar} for more background on $\nabla_A^{Hom}$). Often, 
we will write $\nabla_A$ for $\nabla_A^{Hom}$ if it is clear we are using 
this induced connection.  The connection  
$\nabla_A^{Hom}$ induces a connection on $\kb$, and allows us to calculate 
$\nabla_A g$ for $g\in\Gamma(\Kb)$.  Note that $\nabla_A g$ is not 
necessarily a section of $\Kb$, but a section of $\mathrm{Hom}(E,E)$.  

Of special interest is the trivial bundle $P:=\bar{O}\times K\to\bar{O}$, where $O$ is a bounded open 
subset of $\mathbb{R}^n$ with smooth boundary.  
In this case, the induced bundles $E$, $\kb$, and $\Kb$ are also direct 
products of the appropriate sort.  For example, $E=\bar{O}\times V\to \bar{O}$.  
In this case, the \emph{flat connection} on $P$ is the Ehresmann connection whose 
kernel is the tangent bundle of $\bar{O}$ embedded in the tangent bundle of $\bar{O}\times K$.  
Using parallel transport, one can check that the Koszul connection $\nabla_0$ induced on the product bundle $E$ is 
given by
\begin{equation*}
\nabla_0 \sigma = d \sigma\mbox{, $\sigma\in\Gamma(E)$}
\end{equation*}
where $d \sigma$ is the standard exterior derivative.  Thus we will often use $\nabla_0$ and $d$ interchangeably  
and call them the \emph{flat connection} on $E$.

We are only concerned with connections $\nabla_A$ on $E$ that are 
induced by connections on $P$.  Such connections are called 
\emph{$K$-connections}, and one can show that $\nabla_A$ is a $K$-connection 
if and only if the local connection form is $\mathfrak{k}$-valued.

As a vector bundle, we may equip $\kb$ with a metric.  Any Ad-invariant 
inner product on $\mathfrak{k}$ will induce a Riemannian metric on $\kb$.  In particular, 
we can use the trace inner product $(A,B)=\tr (A^*B)$ to induce a metric on 
$\kb$.  We now view $\kb$ 
as equipped with the metric induced by the trace inner product on 
$\mathfrak{k}$.  Similarly, we equip the bundle $\Mb$ with the trace inner 
product.

For any vector bundle $\xi$ over $M$, we may consider  
the associated vector bundle $\mathrm{Hom}(\Lambda^j(TM),\xi)$
for any $j\ge 1$ and define $\mathrm{Hom}(\Lambda^0(TM),\xi):=\xi$.
We call sections of $\mathrm{Hom}(\Lambda^j(TM),\xi)$ \emph{$\xi$-valued $j$-forms}, and 
generally call them \emph{vector-valued forms}. We denote these sections by $\Omega^j(\xi)$. 
Any connection on $\xi$ induces a connection on $\mathrm{Hom}(\Lambda^j(TM),\xi)$ 
that involves the Levi-Civita connection on $M$.  See \cite{Dar} for more 
about these forms and the induced connections. 

There are certain operations we will like to define on forms.  Given any 
$\alpha\in\Omega^1(\kb)$ and $\phi\in\Gamma(\kb)$, we define the $1$-form 
$[\alpha,\phi]$ by
\begin{equation*}
[\alpha,\phi](X) = [\alpha(X),\phi],
\end{equation*}
for any $X\in TM$.  Also, given any $\alpha,\beta\in\Omega^1(\kb)$, we define the product 
$[\alpha\cdot\beta]\in \Gamma(\kb)$ in the following way:  Suppose 
locally $\alpha = \sum_i \alpha_i dx_i$, and $\beta = \sum_i \beta_i dx_i$, 
and the associated metric tensor is $g_{ij}$.  Then, we set
\begin{equation}\label{dotproduct}
[\alpha\cdot\beta] = \sum_{i,j} [\alpha_i,\beta_j]g^{ij},
\end{equation}
noting that $<dx_i,dx_j>=g^{ij}$ where the matrix $(g^{ij})$ is inverse to 
$(g_{ij})$ .  
One can verify that this globally defines $[\alpha\cdot\beta]$ as a section of $\kb$.

We will often be looking at the difference between two $K$-connections, 
and the following will be useful in looking at such differences. If 
$\nabla_{A_1}$ and $\nabla_{A_2}$ are two
$K$-connections, using the local characterization of $K$-connections, one 
can show that the difference 
$\nabla_{A_1}-\nabla_{A_2}$ is a $\kb$-valued $1$-form.  Furthermore, 
if we set $\alpha := \nabla_{A_1}-\nabla_{A_2}$, we have for any 
$\phi\in\Gamma(\kb)$
\begin{equation}\label{dprop}
(\nabla_{A_1}^{Hom}-\nabla_{A_2}^{Hom})(\phi)=[\alpha,\phi].
\end{equation}
Similarly, if $\beta\in\Omega^1(\kb)$, one can show that
\begin{equation}\label{d*prop}
((\nabla_{A_1}^{Hom})^*-(\nabla_{A_2}^{Hom})^*)(\beta) = -[\alpha\cdot\beta].
\end{equation}
The previous two equations are ubiquitous in what follows.  On sections we 
have $d_A=\nabla_A$ and on $1$-forms we have $d_A^*=(\nabla_A)^*$. We will 
use both notations interchangably on these respective domains.
The curvature $R_A$ of a $K$-connection $\nabla_A$ is a 
$\kb$-valued $2$-form.

Using \eqref{dprop}, one can show that a $K$-connection $\nabla_A^{Hom}$ is 
compatible with the metric on $\Mb$ induced by the trace inner product. Furthermore, a 
$K$-connection $\nabla_A$ on $E$ and the Levi-Civita connection on $M$ 
induce a connection on 
$\mathrm{Hom}(\Lambda^j(TM),\kb)$ that is compatible with the induced metric on 
$\mathrm{Hom}(\Lambda^j(TM),\kb)$.
%
%
%
%
\subsection{Dirichlet Boundary Conditions and Sobolev Spaces}\label{1.2}
We define Sobolev spaces of sections of vector bundles as Palais does in 
\cite{Palais}.  Using the notation of \cite{Palais}, the space $L^p_k(\xi)$ is the 
space of sections of $\xi$ with $k$ Sobolev derivatives under the $L^p$ norm, and $L^p_k(\xi)^0$ is the 
completion of $C^\infty_c(\xi|_{\mathrm{int}(M)})$ in the $L^p_k$ norm.  As usual, we define $H^k(\xi) := L^2_k(\xi)$, where the latter 
notation is what \cite{Palais} uses.  Also converting from Palais's notation, 
we set $H^k_0(\xi) := L^2_k(\xi)^0$.

Marini in \cite{Ma} defines Dirichlet boundary conditions on Sobolev 
spaces, which we will denote $H^k_{Dir}(\mathrm{Hom}(\Lambda^j(TM),\xi))$ for appropriate vector 
bundles $\xi$.  Specifically, 
we define the \emph{Dirichlet Sobolev space $H^k_{Dir}(\mathrm{Hom}(\Lambda^j(TM),\xi))$} for $k\ge 1$ as
\begin{eqnarray*}
H^k_{Dir}(\mathrm{Hom}(\Lambda^j(TM),\xi))&:=& \{\alpha\in H^k(\mathrm{Hom}(\Lambda^j(TM),\xi)): \iota^*(\alpha)=0, \\
& &\mbox{ where $\iota:\partial M\to M$ is the inclusion}\}.
\end{eqnarray*}
Since $k\ge 1$, $\alpha|_{\partial M}$ is defined in the trace sense, so $\iota^*(\alpha)$ is 
defined almost everywhere.  Note that $\iota^*(\alpha)=0$ is equivalent to saying that $\alpha$ vanishes 
on wedges of vectors $X_1\wedge\ldots\wedge X_j$, where all $X_i$ are tangent to $\partial M$.  
For a $0$-form $\sigma$ (i.e. a section $\sigma$), $\iota^*(\sigma)=0$ if and only if $\sigma|_{\partial M}=0$.
Hence, we see that
\begin{equation}
H^k_{Dir}(\xi) = H^1_0(\xi)\cap H^k(\xi)\mbox{, $k\ge 1$}\label{pB2}.
\end{equation}

In what follows we use $k>n/2 +1$ where $k$ is an integer so we can use the multiplication theorem of 
Sobolev spaces for $H^{k-1}$ (see Corollary 9.7 in \cite{Palais}).

Since we will be using Dirichlet boundary conditions, we need a fixed smooth 
connection $\nabla_{A_0}$.  Set 
\begin{equation*}
\Ck := \nabla_{A_0} + H^k_{Dir}(\mathrm{Hom}(TM,\kb)).
\end{equation*}
Note that all the connections in $\Ck$ will be equal to $\nabla_{A_0}$ in tangential 
directions on the boundary.  Also $\Ck$ is an affine space and is therefore a $C^\infty$-Hilbert manifold.  We will call any connection $\nabla_{A}$ 
\emph{$C^\infty$-smooth} (resp. \emph{$H^k$-smooth, $L^p$-smooth}) if $\nabla_A-\nabla_{A_0}\in C^\infty(\mathrm{Hom}(TM,\kb))$ (resp.
$\in H^k(\mathrm{Hom}(TM,\kb))$, $L^p(\mathrm{Hom}(TM,\kb))$.  Note that if $\nabla_{A}$ is $C^\infty$-smooth, then it is a Koszul connection in the usual Riemannian geometry sense.

The sections of $\Kb$ are the \emph{gauge transformations}.  The Sobolev regularity and boundary conditions we 
will need is set in the following definition:

\begin{definition}
Let $\nabla_{A_0}$ be a fixed smooth $K$-connection on $E$.  Let $g\in H^{k+1}(\Kb)$, with 
$g|_{\partial M}\equiv e$, where $e$ is the identity on $K$.  Then we say that $g\in\Gk$, 
and call $\Gk$ the \emph{(Dirichlet) gauge group}.
\end{definition}

The Sobolev space $H^{k+1}(\Kb)$ is defined as in \cite{MV} as the 
completion of smooth sections of $\Kb$ in the norm $H^{k+1}(\Mb)$.  This 
completion without 
the boundary conditions we will call $\mathcal{G}^{k+1}$, as it is called 
in \cite{MV}.

\begin{proposition}\label{blgroup}
For $k>n/2+1$, the group $\Gk$ is a Hilbert Lie group whose Lie algebra is identifiable with $H^{k+1}_{Dir}(\kb)$.  

\end{proposition}
\begin{proof}
As proven in \cite{MV}, $\mathcal{G}^{k+1}$ is Hilbert Lie group.  Since $\Gk$ is a closed topological subgroup of $\mathcal{G}^{k+1}$, it carries a topology.  To give $\Gk$ Hilbert coordinates in a neighborhood of the identity, we show that the exponential map takes $H^{k+1}_{Dir}(\kb)$ into $\Gk$ and is a local 
homeomorphism at $0$.  In \cite{MV} it is shown that $\exp$ is a $C^\infty$ local diffeomorphism $\exp: H^{k+1}(\kb)\to\mathcal{G}^{k+1}$ without boundary conditions.  Hence, 
we need only show that $\exp$ maps $H^{k+1}_{Dir}(\kb)$ into $\Gk$, and for a neighborhood $U$ of the 
identity in $\mathcal{G}^{k+1}$, $\exp^{-1}\equiv\log$ maps $U$ into $H^{k+1}_{Dir}(\kb)$.

To prove the first assertion, let $f\in H^{k+1}_{Dir}(\kb)$.  Since $k+1>n/2 + 2$ we have 
$f\in C^2(\kb)$ with $f|_{\partial M}\equiv 0$.  Then if $g:=\exp(f)$, we have that
$g\in C^2(\Kb)\subseteq C^2(\Mb)$, and 
$g|_{\partial M}\equiv e$, where $e$ is the identity element of $K$, proving $g\in\Gk$.

To prove the second assertion, since $\exp:H^{k+1}(\kb)\to\mathcal{G}^{k+1}$ 
is a local 
diffeomorphism between the spaces without boundary conditions, we need only 
show that for small $\mu\in H^{k+1}(\kb)$, if $g:=\exp(\mu)\in\Gk$, 
then $\mu\in H^{k+1}_{Dir}(\kb)$.  Note that $\sup |\mu|\le C\|\mu\|_
{H^{k+1}}$. 
So for small enough $\|\mu\|_{H^{k+1}}$, we can use the fact that the ``pointwise'' 
map $\exp:\mathfrak{k}\to
K$ is local diffeomorphism at $0$ to say that since $g|_{\partial M}\equiv e$, 
we have $\mu|_{\partial M}\equiv 0$.  Since $\mu\in H^{k+1}(\kb)$, this 
vanishing on the boundary implies that $\mu\in H^1_0(\kb)\cap H^{k+1}(\kb)=
H^{k+1}_{Dir}(\kb)$.

To give $\Gk$ an atlas, we transport these coordinates via right translation as is done in \cite{MV}.
\end{proof}

$\Gk$ acts on $\Ck$ on the right in the following way.  Suppose that we have a $1$-form $\eta\in H^k_{Dir}(\mathrm{Hom}(TM,\kb))$.  Then the action is
\begin{equation}
(\nabla_{A_0} + \eta)\cdot g = \nabla_{A_0} + (g^{-1}\nabla^{Hom}_{A_0}g + 
\mathrm{Ad }(g^{-1})\eta).\label{gaction}
\end{equation}
By the same reasoning as found in \cite{MV}, this action is smooth.
Note that for $(\nabla_{A_0} + \eta)\cdot g$ to remain in $\Ck$, we need to have $g^{-1}\nabla^{Hom}_{A_0}g $ satisfy Dirichlet boundary conditions.  The following proposition shows that this is the case.

\begin{proposition}\label{gaugeEZ}
Let $k>n/2+1$ and suppose $g\in\Gk$.  Then we have $g^{-1}\nabla_{A_0}g\in H^k_{Dir}(\mathrm{Hom}(TM,\kb))$, 
\end{proposition} 
\begin{proof}
Let $g\in\Gk$.  Since $g\in H^{k+1}(\Kb)$, there exist smooth $g_m\in H^{k+1}(\Kb)$ such that $g_m\to g$ in 
$H^{k+1}(\Mb)$.  It is shown in \cite{MV} that inversion is continuous on $\mathcal{G}^{k+1}$.
Hence, $(g_m)^{-1}\to g^{-1}$ in $H^{k+1}(\Mb)$.  Since $\nabla_{A_0}$ is a smooth $K$-connection, we see that 
$g_m^{-1}\nabla_{A_0}^{Hom} g_m\in\Omega^1(\kb)$, and by the multiplication theorem,
\begin{eqnarray*}
||g^{-1}\nabla_{A_0}^{Hom} g - g_m^{-1}\nabla_{A_0}^{Hom} g_m||_{H^k} &\le& ||g^{-1}\nabla_{A_0}^{Hom} g - g_m^{-1}\nabla_{A_0}^{Hom} g||_{H^k}+\\ & &||g_m^{-1}\nabla_{A_0}^{Hom} g - g_m^{-1}\nabla_{A_0}^{Hom} g||_{H^k}\\
&\le& C(||g^{-1}-g_m^{-1}||_{H^{k+1}}||\nabla_{A_0}^{Hom}g||_{H^k} + \\
& & ||g_m^{-1}||_{H^{k+1}}||g_m - g||_{H^{k+1}})\to 0.
\end{eqnarray*}
Thus, $g^{-1}\nabla_{A_0}^{Hom} g\in H^k(\mathrm{Hom}(TM,\kb))$.  

We now show that $g^{-1}\nabla_{A_0}^{Hom} g\in H^k_{Dir}(\mathrm{Hom}(TM,\kb))$.  Locally, we have $\nabla_{A_0} = d + A_0$, where $A_0$ is a $C^\infty$-smooth $\mathfrak{k}$-valued $1$-form.  Let $X$ be a tangential direction at a boundary point.
Since $g\equiv e$ on $\partial M$, we have $dg(X)=0$. Also on the boundary, $[A_0(X),g]=[A_0(X),e]=0$, since 
$e$ commutes with everything.  Hence, globally, $\iota^*(\nabla_{A_0}^{Hom} g) \equiv 0$, and thus 
$\iota^*(g^{-1}\nabla_{A_0}^{Hom} g) \equiv 0$.  This proves that $g^{-1}\nabla_{A_0}^{Hom} g\in H^k_{Dir}(\mathrm{Hom}(TM,\kb))$.
\end{proof}

Using Dirichlet boundary conditions gives us a Sobolev-\Poincare inequality.

\begin{proposition}\label{SPprop}
Let $\nabla_{A}$ be a $L^q$-smooth connection on a vector bundle $\xi\to M$ compatible with 
the metric on $\xi$, where $q\ge n$.  Then there exists a $\kappa_p > 0$ such that for any 
$f\in L^p_1(\xi)^0$ with $1<p<n$, we have
\begin{equation}
\|f\|_{L^p}\le\kappa_p\|\nabla_A f\|_{L^p}\label{SP},
\end{equation}
where $\kappa_p$ does not depend on the connection, but does depend on $p$.
\end{proposition}
\begin{proof}
This is done for real-valued functions by a standard \Poincare inequality argument that can be found in, 
for example, \cite{Evans} and \cite{He}.  This shows that there exists a $\kappa_p >0$ such that 
\begin{equation}
\|g\|_{L^p}\le\kappa_p\|dg\|_{L^p},\label{HSP}
\end{equation}
for any $g\in L^p_1(M\times\mathbb{R})^0$ where $M\times\mathbb{R}\to M$ is the trivial vector bundle.  (The references 
above prove \eqref{HSP} for real-valued functions $g$.  But real-valued functions on $M$ and sections of $M\times\mathbb{R}$ 
are the same.)

Now let $f\in C^1_0(\xi|_{\mathring{M}})$, where $\mathring{M}$ is the interior of $M$.  Then the function 
$|f|$ is globally Lipschitz, so by Lemma 2.8 in \cite{He} we have 
$|f|\in L^p_1(M\times\mathbb{R})$.  Since $|f|$ is continuous and $0$ 
on $\partial M$, Theorem 5.5.2 in \cite{Evans} tells us that $|f|\in L^p_1(M
\times\mathbb{R})^0$.  Hence, \eqref{HSP} yields
\begin{eqnarray*}
\|f\|_{L^p}\le\kappa_p\|d|f|\|_{L^p}\le\kappa_p\|\nabla_A f\|_{L^p}.
\end{eqnarray*}
The second inequality is \emph{Kato's inequality}.  This inequality only requires that $\nabla_A$ is 
compatible with the metric.  For a proof of this inequality, see \cite{Par}.
Since $C^1_0(\xi|_{\mathring{M}})$ is dense in $L^p_1(\xi)^0$, the proposition 
has been proven. (The condition $q\ge n$ ensures that $\nabla_A f \in L^p$).
\end{proof}

The Sobolev-\Poincare inequality immediately tells us

\begin{corollary}\label{gfree}
For $k > n/2 + 1$, the action of $\Gk$ on $\Ck$ is free.
\end{corollary}
\begin{proof}
Suppose $\nabla_{A_0} + \eta\in\Ck$
, $g\in\Gk$, and 
$(\nabla_{A_0} +\eta)\cdot g = \nabla_A + \eta$.  By (\ref{gaction}), this means 
that
\begin{equation*}
(\nabla^{Hom}_{A_0} +\eta)g = \nabla^{Hom}_{A_0} g + [\eta,g] = 0.
\end{equation*}
By (\ref{SP}), we have
\begin{equation*}
\|g-e\|_{L^2} \le \kappa_2\|(\nabla^{Hom}_{A_0} +\eta)g\|_{L^2} = 0
\end{equation*}
Since $g$ is continuous, the above shows $g\equiv e$ and thus the corollary is proven.
\end{proof}

The freeness of this action allows us to directly use $\Ck$ and $\Gk$, instead of 
so-called generic (or irreducible) connections and modified gauge groups as found in \cite{MV}, \cite{NR}, and \cite{Par}.

Given a $K$-connection $\nabla_{A}\in\Ck$, we can define the Laplacian
\begin{equation*} 
\Delta_A=\nabla_A^*\nabla_A = d^*_A d_A:H^{m+1}_{Dir}(\kb)\to H^{m-1}(\kb)\mbox{ for $1\le m\le k$.}
\end{equation*}   
The regularity is correct by the following:  Since $\nabla_{A_0}$ is 
a smooth connection, clearly $\Delta_{A_0}$ is a bounded map from $H^{m+1}$ into 
$H^{m-1}$.  Suppose $h=\nabla_A-\nabla_{A_0}\in H^k_{Dir}(\kb)$.  Recalling 
equations \eqref{dprop} and \eqref{d*prop} and the comment following them, 
for $f\in H^{m+1}$ we have
\begin{equation}\label{localLap}
\Delta_{A} f = \Delta_{A_0}f + [d_{A_0}^*h,f] - [h\cdot[h,f]] - 2[h\cdot 
d_{A_0}f].
\end{equation}
So, we have (allowing $||\cdot||_i$ to denote the $H^i$ norm)
\begin{eqnarray*}
||\Delta_A f||_{m-1} &\le& ||\Delta_{A_0}f||_{{m-1}}+||[d_{A_0}^*h,f]||_
{{m-1}} \\
& &+ ||[h\cdot[h,f]]||_{m-1} + 2||[h\cdot d_{A_0}f]||_{m-1} \\
&\le& ||\Delta_{A_0}f||_{{m-1}}+C(||d_{A_0}^*h||_{k-1}||f||_{{m-1}} \\
& &+ ||h||_k||[h,f]||_{m-1} + 2||h||_k||d_{A_0}f]||_{m-1}) \\
&\le& ||\Delta_{A_0}f||_{{m-1}}+C(||h||_{k}||f||_{{m-1}} \\
& &+ ||h||_k^2||f||_{m-1} + 2||h||_k||f||_{m}) < \infty,
\end{eqnarray*}
where we used the fact that $H^{m-1}$ is a $H^{k-1}$ module, which is the case 
since $k-1\ge m-1$ and $k-1 > n/2$ (this is a critical point where we need $k > n/2+1$).  Thus, $\Delta_A$ is bounded from $H^{m+1}_
{Dir}(\kb)$ to $H^{m-1}(\kb)$. Furthermore, we have
\begin{proposition}\label{Greenexist}
Let $\nabla_{A}\in\Ck$ for an integer $k > n/2 + 1$, and suppose $1\le m\le k$.  Then the mapping 
$\Delta_A:H^{m+1}_{Dir}(\kb)\to H^{m-1}(\kb)$ is an isomorphism.  Furthermore, 
if $f\in H^2_{Dir}(\kb)$ and $\Delta_A f\in H^{m-1}(\kb)$, then $f\in 
H^{m+1}_{Dir}(\kb)$ and 
\begin{equation}
\|f\|_{H^{m+1}}\le C(\|\Delta_A f\|_{H^{m-1}} + \|f\|_{H^0}).
\end{equation}
\end{proposition}
We set $G_A:=(\Delta_A)^{-1}$ and call it a \emph{Green operator}.  

The proof of the existence of Green operators follows from variations of standard argument for weak solutions 
to elliptic equations that can be found in \cite{Evans} and \cite{GT}. We omit the details although they can 
be found in \cite{GrycThesis}. Nominally, the proof in \cite{GrycThesis} is for $n=3$.  But replacing ``$3$'' with ``$n$'' in the proof gives the general result.  

We emphasize here that 
\emph{every} connection $\nabla_A\in\Ck$ has an associated Green operator.  We need not restrict our space of connections in 
this with-boundary case since we are imposing Dirichlet boundary conditions.
%
%
%
%
\subsection{The Quotient $\Ck\to\Ck/\Gk$}\label{1.3}
We are now in a position to consider the structure of the quotient $\Ck\to\Ck/\Gk$.
\begin{proposition}\label{bigpB}
Let $k>n/2+1$ for an integer $k$.  The quotient space $\Ck/\Gk$ is a $C^\infty$ Hilbert manifold, and $\pi: \Ck\to 
\Ck/\Gk$ is a principal bundle with structure group $\Gk$.
\end{proposition}
The proof is a straightforward adaptation of the no-boundary case of Mitter and Viallet in \cite{MV}.
Others have proven similar statements in more specific situations (see \cite{AHS}, \cite{NR}, and \cite{Par}).
Since the techniques for proving Proposition \ref{bigpB} are very little changed from those employed by the above authors, we omit the proof.  However, for those interested in the proof, it is in \cite{GrycThesis}.  Again, the proof in \cite{GrycThesis} is nominally for $n=3$, but an examination of the proof shows that the dimension is not mentioned and the only relevent contribution of the dimension is that we have $k>n/2+1$ as we have here.  Thus, the proof works for the $n\ge 2$ as well.
We again note that with our Dirichlet boundary conditions on the connections and the gauge group, we need not restrict the 
space of connections to generic connections nor restrict the gauge group further.  This situation is unlike the no-boundary situations as found in \cite{AHS}, \cite{NR}, and \cite{Par}.  Also, the presence of the Sobolev-\Poincare inequality of Proposition \ref{SPprop} leads to some simplifications of the arguments.
%
%
%
%
\section{The Coulomb Connection and Its Holonomy}\label{Coulomb}
Since the bundle $\Ck\to \Ck/\Gk$ is a principal bundle, 
we can consider holonomy groups of a fixed connection upon it.  The connection we will consider is the \emph{Coulomb connection} whose connection form at $\nabla_A$ is defined as $\displaystyle G_A d^*_A$.  The corresponding horizontal 
subspace at $\nabla_{A}$ we will denote as $H_A$.  Note that since $\Ck = \nabla_{A_0} + H^k_{Dir}(\mathrm{Hom}(TM,\kb))$, the 
tangent space at $\nabla_A\in\Ck$ is \\
$H^k_{Dir}(\mathrm{Hom}(TM,\kb))$.  One can verify that $H_A$ is
\begin{equation*}
H_A = \{\alpha\in H^k_{Dir}(\mathrm{Hom}(TM,\kb)): d_A^*\alpha = 0\}.
\end{equation*}
This connection is natural in the sense that $H_A$ is the $L^2$ orthogonal 
complement to the vertical vectors at $\nabla_{A}$.  Indeed, one can show that given $\gamma\in 
\Lie(\Gk)=H^{k+1}_{Dir}(\mathrm{Hom}(TM,\kb))$, the fundamental vector field associated to $\gamma$ is  
$d_A \gamma$. Hence the vertical vectors are those vectors of the form $d_A\gamma$ for some 
$\gamma\in\Lie(\Gk)$ (see \cite{Gross} or \cite{NR}). 
By the same reasoning as the proof of Lemma 7.1 in \cite{NR}, the Coulomb connection is indeed 
a connection on $\Ck\to\Ck/\Gk$. 

The principal tool of our study of the holonomy group of the Coulomb connection is the image of 
the curvature form of the Coulomb connection.  Let $\mathcal{R}_A$ be the 
curvature form of the Coulomb connection at $\nabla_A$.  By the same calculation 
in the proof of Lemma 7.2 in \cite{NR}, we have 
\begin{equation}\label{curveq}
\mathcal{R}_A(\alpha,\beta)=-2G_A([\alpha\cdot\beta]),\mbox{ for } \alpha,\beta\in H_A.
\end{equation}
\subsection{Coordinates at the Boundary of $\partial M$}\label{2.1}
In this investigation, certain types of coordinates at the boundary are useful.
Consider the following system of coordinates at the boundary that satisfy the following:
\begin{enumerate}
\item[A1.] $\partial/\partial x_n$ is orthogonal to $\partial/\partial x_1,\ldots,\partial/\partial x_{n-1}$ on the boundary.\label{A1}
\item[A2.] $\partial/\partial x_n$ has norm 1 everywhere.\label{A2}
\item[A3.] $\partial/\partial x_n$ is the inward pointing unit normal vector on the 
boundary.\label{A3}
\end{enumerate}
We describe such a coordinate system as \emph{Type A}.  They have been constructed in \cite{Milnor1} and \cite{RS}.  We also will use a similar coordinate system such that
\begin{enumerate}
\item[B1.] $\partial /\partial y_n$ is orthogonal to $\partial /\partial y_1,\ldots\partial /\partial y_{n-1}$ 
everywhere.
\item[B2.] $\partial /\partial y_n$ is a positive (perhaps nonconstant) multiple of the 
inward pointing unit normal vector on the boundary.
\end{enumerate}
We call such coordinates \emph{Type B}.  A detailed contruction of these coordinates can be found in \cite{GrycThesis}.  In what follows, we will use $\{x_1,\ldots,x_n\}$ to denote Type A coordinates and $G=(g_{ij})$ to 
denote the associated metric tensor. For Type B coordinates, 
we use $\{y_1,\ldots,y_n\}$ and $\tilde{H}=(h_{ij})$, respectively (we use $\tilde{H}$ to distinguish this matrix from the mean curvature of the boundary, which we will denote $H_{}$).

\subsection{Mean Curvature of the Boundary and Coordinates}\label{2.2}
We will see that the mean curvature $H_{}$ of the immersion $\iota:\partial M\to M$ will come into play in our investigation of 
the image of the curvature form.  (For the definition of mean curvature, see, for example, \cite{doCarmo}).)  It will be useful to have a characterization of $H_{}$ in our Type A and Type B coordinates.
\begin{proposition}\label{taux}
Consider the mean curvature $H_{}$ of the immersion $\iota:\partial M\to M$.  Let $\{x_1,\ldots,x_n\}$ be Type A coordinates with 
metric tensor $G=(g_{ij})$.  Then $H_{}$ can be written as 
\begin{equation}
H_{}(x_1,\ldots,x_{n-1}) = \frac{1}{n-1}\frac{\partial a}{\partial x_n}(x_1,\ldots,x_{n-1},0)\cdot\frac{1}{a(x_1,\ldots,x_{n-1},0)},
\end{equation}
where $a:=\sqrt{\det(g_{ij})}$.\footnote{For those also reading \cite{GrycThesis} the $H_{}$ defined here is $1/2$ times the 
$\tau$ defined in \cite{GrycThesis}.}
\end{proposition}
\begin{proof}
Let $\Gamma^m_{ij}$ be the Christoffel symbols corresponding to $\{x_1, \ldots, x_n\}$.  Since the connection we are considering is the Levi-Civita connection, we have
\begin{equation*}
\Gamma^m_{ij}=\frac{1}{2}\sum_{k=1}^n \left(\frac{\partial}{\partial x_i}(g_{jk})+\frac{\partial}{\partial x_j}(g_{ik})- \frac{\partial}{\partial x_k}(g_{ij})\right)g^{km}.
\end{equation*}
(see, for example, \cite{doCarmo}).  By our choice of coordinate system, we have
\begin{equation*} 
g_{1n}=\ldots=g_{(n-1)n}=g^{1n}=\ldots=g^{(n-1)n}=0
\end{equation*}
on the boundary and $\frac{\partial}{\partial x_1},\ldots,\frac{\partial}{\partial x_{n-1}}$ are tangent to the boundary.  So on the boundary we have for $i<n$
\begin{eqnarray}
\Gamma^i_{in}=\frac{1}{2}\sum_{k=1}^{n-1}\left(\frac{\partial}{\partial x_n}(g_{ik})g^{ki}\right)\label{tau1}.
\end{eqnarray}
Using a Laplace expansion on the bottom row of $G$, we also have \begin{equation}\label{laplace}
\det(G) = \sum_{k=1}^n (-1)^{n+k} g_{nk}\cdot\det(G(n|k)).\end{equation}
Above and in what follows, $G(i|j)$ is the $(n-1)$-by-$(n-1)$ matrix obtained from $G$ by deleting the $i^\mathrm{th}$ row and $j^\mathrm{th}$ column.  On the boundary, for $k<n$ we have $g_{nk}=0$ and $\det(G(n|k))=0$, since the last column of $G(n|k)$ is all zeros. Also, $g_{nn}\equiv 1$ everywhere by our choice of Type A coordinagtes.  Using these preceeding facts, we apply the product rule to \eqref{laplace} to get
\begin{equation}\label{tau4}
\frac{\partial}{\partial x_n}\det(G) = \frac{\partial}{\partial x_n}\det(G(n|n))
\mbox{  on $\partial M$.}
\end{equation}
Define a set of permuations $S^{i,j}_{n-1}$ as
\begin{equation*}
S^{i,j}_{n-1} = \{\sigma\in S_{n-1}: \sigma(i)=j\}.
\end{equation*}
Then $S^{i,j}_{n-1}$ is isomorphic to $S_{n-2}$, and if $\sigma\in S^{i,j}_{n-1}$ corresponds to $\tilde{\sigma}\in S_{n-2}$, then one can show that
\begin{equation}\label{tau5}
\mathrm{sgn}(\sigma)=(-1)^{(i+j)}\mathrm{sgn}(\tilde{\sigma}).
\end{equation}
Indeed, one can prove \eqref{tau5} by considering the permuation matrix of $\sigma$ (the determinant of which is $\mathrm{sgn}(\sigma)$), and using a Laplace expansion down the $i^{th}$ row.  Define the $(n-1)$-by-$(n-1)$ matrix  $C:=G(n|n)$.  Combining \eqref{tau1}) and \eqref{tau4} we have 
\begin{eqnarray*}
H_{} &=& -\frac{1}{n-1}\sum_{i=1}^{n-1}\left(\nabla_{\frac{\partial}{\partial x_i}}(-\frac{\partial}{\partial x_n})\right)_i\\
&=& \frac{1}{n-1}\sum_{i=1}^{n-1} \Gamma^i_{in} = \frac{1}{n-1}\sum_{i=1}^{n-1}\frac{1}{2}\sum_{j=1}^{n-1}\frac{\partial}{\partial x_n}(g_{ij})g^{ji}\\
&=& \frac{1}{2(n-1)}\sum_{i=1}^{n-1}\sum_{j=1}^{n-1}\frac{\partial}{\partial x_n}(g_{ij})c^{ji}\\
&=& \frac{1}{2(n-1)\det(C)}\sum_{i=1}^{n-1}\sum_{j=1}^{n-1}\frac{\partial}{\partial x_n}(g_{ij})(-1)^{(i+j)}\det(C(i|j))\\
&=& \frac{1}{2(n-1)\det(C)}\sum_{i=1}^{n-1}\sum_{j=1}^{n-1}\sum_{\tilde{\sigma}\in S_{n-2}} \frac{\partial}{\partial x_n}(g_{ij})(-1)^{(i+j)}\mathrm{sgn}(\tilde{\sigma})c_{1\tilde{\sigma}(1)}\cdot\ldots\cdot\widehat{c_{ij}}\cdot\\
& & \ldots\cdot c_{(n-1)\tilde{\sigma}(n-1)}. \\
\end{eqnarray*}
Inserting \eqref{tau5} to the above yields
\begin{eqnarray*}
H&=& \frac{1}{2(n-1)\det(C)}\sum_{i=1}^{n-1}\sum_{j=1}^{n-1}\sum_{\sigma\in S^{i,j}_{n-1}}\mathrm{sgn}(\sigma) \frac{\partial}{\partial x_n}(g_{ij})c_{1\sigma(1)}\cdot\ldots\cdot\widehat{c_{ij}}\cdot\\
& & \ldots\cdot c_{(n-1)\sigma(n-1)} \\
&=& \frac{1}{2(n-1)\det(G(n|n))}\sum_{i=1}^{n-1}\sum_{\sigma\in S_{n-1}} \frac{\partial}{\partial x_n}(g_{i\sigma(i)})g_{1\sigma(1)}\cdot\ldots\cdot\widehat{g_{i\sigma(i)}}\cdot\\
& & \ldots\cdot g_{(n-1)\sigma(n-1)} \\
&=& \frac{1}{2(n-1)\det(G(n|n))}\sum_{\sigma\in S_{n-1}} \frac{\partial}{\partial x_n}(g_{1\sigma(1)}\cdot\ldots\cdot g_{(n-1)\sigma(n-1)}) \\
&=& \frac{1}{2(n-1)\det(G(n|n))}\frac{\partial}{\partial x_n}(\det(G(n|n))) = \frac{1}{2(n-1)\det(G)}\frac{\partial}{\partial x_n}(\det(G)) \\
&=& \frac{1}{(n-1)}\frac{\partial a}{\partial x_n}\frac{1}{a},
\end{eqnarray*}
as desired.
\end{proof}

We can also write $H_{}$ in terms of Type B coordinates:
\begin{proposition}\label{tauy}
Consider the mean curvature $H_{}$ of the immersion $\iota:\partial M\to M$.  Consider Type B coordinates $\{y_1,\ldots,y_n\}$ with metric tensor $\tilde{H}=(h_{ij})$.  Then $H_{}$ can be written as
\begin{equation*}
H_{} = \frac{1}{n-1}\left({\sqrt{h_{nn}}}\cdot d\left(\frac{1}{\sqrt{h_{nn}}}\right)(\nu) + \frac{db(\nu)}{b}\right),
\end{equation*}
where $b:=\sqrt{\det(h_{ij})}$ and $\nu$ is the unit inward pointing normal vector.
\end{proposition}
\begin{proof}
One can use a coordinate change between Type A and Type B coordinates (see the $n=3$ case in \cite{GrycThesis}).  Or, one can proceed similarly as the proof of Proposition \ref{taux}.  Indeed, we have a similar chain of equations from the proceeding proof.  Below we write a shortened list of equations where the steps that are the same in the Type A case are omitted:
\begin{eqnarray*}
H_{} &=& -\frac{1}{n-1}\sum_{i=1}^{n-1}\left(\nabla_{\frac{\partial}{\partial y_i}}(-\frac{1}{\sqrt{h_{nn}}}\frac{\partial}{\partial y_n})\right)_i\\
&=& \frac{1}{\sqrt{h_{nn}}(n-1)}\sum_{i=1}^{n-1}\frac{1}{2}\sum_{j=1}^{n-1}\frac{\partial}{\partial y_n}(h_{ij})h^{ji}\\
&=& \frac{\sqrt{h_{nn}}}{2(n-1)\det(\tilde{H})}\sum_{i=1}^{n-1}\sum_{j=1}^{n-1}\frac{\partial}{\partial y_n}(h_{ij})C(i|j)\\
&=& \frac{\sqrt{h_{nn}}}{2(n-1)\det(\tilde{H})}\frac{\partial}{\partial y_n}(\det(\tilde{H}(n|n)))\\
&=& \frac{\sqrt{h_{nn}}}{2(n-1)\det(\tilde{H})}\frac{\partial}{\partial y_n}\left(\frac{\det(\tilde{H})}{h_{nn}}\right)\\
&=& \frac{h_{nn}}{2(n-1)\det(\tilde{H})}d\left(\frac{\det(\tilde{H})}{h_{nn}}\right)(\nu)\\
&=& \frac{1}{n-1}\left({\sqrt{h_{nn}}}\cdot d\left(\frac{1}{\sqrt{h_{nn}}}\right)(\nu) + \frac{db(\nu)}{b}\right),
\end{eqnarray*}
as desired.  Twice above (once in the beginning and once at the end), we used the fact that $\frac{1}{\sqrt{h_{nn}}}\frac{\partial}{\partial y_n}=\nu$.
\end{proof}
%
%
%
%
\subsection{The Image of $\mathcal{R}_A$}\label{2.3}
We now use mean curvature of $\partial M$ in the following lemma, which relates $H$ and the 
image of the curvature form.
\begin{lemma}\label{bct}
Suppose $M$ is a $n$ dimensional manifold with boundary, $k>n/2+1$ where $k$ is an integer, $\alpha,\beta\in H_A$ and $\nabla_{A}\in\Ck$.  Then 
\begin{equation}\label{thecondition1}
d_A[\alpha\cdot\beta](\nu)= -2(n-1)H_{}[\alpha\cdot\beta]\mbox{    on $\partial M$,}
\end{equation}
where $\nu$ is the unit inward pointing vector field and $H$ is the mean curvature of the boundary.
\end{lemma}

\noindent Since $k-n/2>1$, note that $[\alpha\cdot\beta]$ is $C^1$, and thus $d_A[\alpha\cdot\beta]$ 
is continuous.  Hence, the above equality is true not just in the trace sense, but as 
an equality of two continuous functions.
\begin{proof}
We will use Type A coordinates $\{x_1,\ldots,x_n\}$, and assume 
that the vector bundle $\k_P$ is also trivialized in this neighborhood.  
Recall that the metric tensor in this coordinate system has the feature that 
$g_{in}=\delta_{in}$ on the boundary, and $g_{nn}=1$ everywhere.  Thus, 
$g^{in}=\delta_{in}$ also on the boundary.  Also, $\frac{\partial}
{\partial x_n}$ is the inward pointing 
unit normal vector on the boundary.  
Take $\alpha,\beta$ as above and define 
$\alpha_i$ and $\beta_i$ so that $\alpha = \sum_{i=1}^n\alpha_i dx_i$ and 
$\beta = \sum_{i=1}^n\beta_i dx_i$.  Since we are assuming $\k_P$ has a fixed 
trivialization in our neighborhood, we can view the $\alpha_i$ and $\beta_i$ 
as $\k$-valued functions.
Also, since $\frac{\partial}{\partial x_n}$ 
is the inward pointing unit normal vector and $\alpha,\beta$ satisfy Dirichlet 
boundary conditions, we have 
\begin{equation}\label{bc1}
\alpha_i=\beta_i=0\mbox{   for $i<n$ on $\partial M$.}
\end{equation}
Let $d$ be the flat connection with respect to our fixed trivialization of 
$\k_P$ and define a $\k$-valued $1$-form $A$ 
so that $d_A = d + A$.  Define $A_i$ so that $A=\sum_{i=1}^n A_idx_i$. On 
this coordinate patch, we have
\begin{eqnarray}\label{bc2}
[\alpha\cdot\beta] &=& \sum_{j,k=1}^n [\alpha_j,\beta_k](dx_j\cdot dx_k)
= \sum_{j,k=1}^n [\alpha_j,\beta_k]g^{jk}.
\end{eqnarray}
Taking the derivative $d_A$ yields
\begin{eqnarray*}
d_A([\alpha\cdot\beta])&=& \sum_{j,k=1}^n d_A([\alpha_j,\beta_k]g^{jk})\\
&=& \sum_{j,k=1}^n [d_A(\alpha_j),\beta_k]g^{jk} + [\alpha_j,d_A(\beta_k)]g^{jk} + [\alpha_j,\beta_k]d(g^{jk})
\end{eqnarray*}
By the properties of Type A coordinates and by \eqref{bc1}, on $\partial M$ the above reduces to
\begin{eqnarray}\label{bc2.5}
d_A([\alpha\cdot\beta])|_{\partial M} = [d_A\alpha_n,\beta_n]+
[\alpha_n,d_A\beta_n] + [\alpha_n,\beta_n]d(g^{nn}).
\end{eqnarray}
Using the adjoint matrix, we see that
\begin{equation*}
g^{nn}=\frac{\det(G(n|n))}{\det(G)}.
\end{equation*}
Combining the fact that $\det(G(n|n))=\det(G)$ on $\partial M$ and \eqref{tau4} yields
\begin{eqnarray}\label{bc2.75}
\frac{\partial g^{nn}}{\partial x_n} &=& \frac{\left(\frac{\partial}{\partial x_n}\det(G(n|n)) \right)\det(G)-
\det(G(n|n))(\frac{\partial}{\partial x_n}\det(G))}{\det(G)^2}\\
&=& 0 \mbox{ on $\partial M$.}
\end{eqnarray}
Hence, we have
\begin{equation}\label{bc3}
d_A([\alpha\cdot\beta])|_{\partial M}(\frac{\partial}{\partial x_n}) = [d_A\alpha_n(\frac{\partial}{\partial x_n}),\beta_n]+
[\alpha_n,d_A\beta_n(\frac{\partial}{\partial x_n})].
\end{equation}

We will leave $d_A([\alpha\cdot\beta])|_{\partial M}$ for the moment 
and investigate what $d^{*}_A\alpha = d^{*}_A\beta = 0$ 
means in our coordinate system.  We will calculate $d^{*}$ by using the 
Hodge star operator.  Set $a:=\sqrt{\det(G)}=1/\sqrt{\det(G^{-1})}$. Using the boundary properties of the $\alpha_i$'s, $\beta_i$'s and $g^{ij}$'s, we have
\begin{eqnarray*}
-d^* \alpha &=& *d*\left(\sum_{i=1}^n\alpha_i dx_i\right) \\
&=& *d\left(a\sum_{i,j=1}^n(-1)^{j+1}\alpha_i g^{ji} dx_1\wedge\ldots\wedge\widehat{dx_j}\wedge\ldots\wedge dx_n \right)\\
&=& *(da)\left(\sum_{i,j=1}^n(-1)^{j+1}\alpha_i g^{ji} dx_1\wedge\ldots\wedge\widehat{dx_j}\wedge\ldots\wedge dx_n \right) + \\
& & \left(\sum_{i,j=1}^n \frac{\partial}{\partial x_i}(\alpha_ig^{ji})\right)*(adx_1\wedge\ldots\wedge dx_n) \\
&=& *[(da)(-1)^{n+1}\alpha_ng^{nn}dx_1\wedge\ldots\wedge dx_{n-1}] + \frac{\partial}{\partial x_n}(\alpha_n)g^{nn} + \\ & & \alpha_n\frac{\partial}{\partial x_n}(g^{nn})\\
&=& \frac{\partial a}{\partial x_n}\frac{1}{a}*(adx_1\wedge\ldots\wedge dx_n) + \frac{\partial}{\partial x_n}(\alpha_n)g^{nn} + \alpha_n\frac{\partial}{\partial x_n}(g^{nn})\\
&=& \frac{\partial a}{\partial x_n}\frac{1}{a} + \frac{\partial \alpha_n}{\partial x_n}g^{nn},
\end{eqnarray*}
where we used \eqref{bc2.75} on the last line.
Since $d_A^*\alpha =d^*\alpha - [A\cdot\alpha]$, we have,
\begin{eqnarray}\label{bc4}
-d^*_A\alpha|_{\partial M}&=& \frac{\partial \alpha_n}{\partial x_n} + 
\frac{1}{a}\frac{\partial a}{\partial x_n}\alpha_n + [A\cdot\alpha]\\
&=& \frac{\partial \alpha_n}{\partial x_n} + 
\frac{1}{a}\frac{\partial a}{\partial x_n}\alpha_n + [A_n, \alpha_n]
\label{bc5},
\end{eqnarray}
where we used (\ref{bc1}) and (\ref{bc2}) (replacing $\beta$ with $A$) in the 
last line.  Of course, an analogous statement holds for $\beta$ replacing $\alpha$.  

We now revisit (\ref{bc3}) and insert (\ref{bc5}):
\begin{eqnarray*}
d_A([\alpha\cdot\beta])|_{\partial M}(\nu) &=& 
[d_A\alpha_n,\beta_n](\partial/\partial x_n)+
[\alpha_n,d_A\beta_n](\partial/\partial x_n) \\
&=& [\frac{\partial\alpha_n}{\partial x_n},\beta_n]+[[A_n,\alpha_n],\beta_n]
\\ & & + [\alpha_n,\frac{\partial\beta_n}{\partial x_n}] + 
[\alpha_n,[A_n,\beta_n]] \\
&=& -[\frac{1}{a}\frac{\partial a}{\partial x_n}\alpha_n + [A_n,\alpha_n]+d_A^*\alpha,
\beta_n]+[[A_n,\alpha_n],\beta_n] \\
& & - [\alpha_n,\frac{1}{a}\frac{\partial a}
{\partial x_n}\beta_n + [A_n,\beta_n]+d_A^*\beta] + [\alpha_n,[A_n,\beta_n]] \\
&=& -\frac{2}{a}\frac{\partial a}{\partial x_n}[\alpha_n,\beta_n]-[d^*_A\alpha,\beta_n]-[\alpha_n, d^*_A\beta] \\
&=& -2(n-1)H_{}[\alpha\cdot\beta]|_{\partial M}-[d^*_A\alpha,\beta(\nu)]-[\alpha(\nu), d^*_A\beta]\\
&=& -2(n-1)H_{}[\alpha\cdot\beta]|_{\partial M}.
\end{eqnarray*}
where we again used (\ref{bc2}) on the second to last line, as well as the fact that $\alpha,\beta\in H_A$.  The lemma is thus proven.
\end{proof}
For future reference, we rewrite the second to last equation above
\begin{equation}
d_A([\alpha\cdot\beta])|_{\partial M}(\nu) = -2(n-1)H_{}[\alpha\cdot\beta]|_{\partial M}-[d^*_A\alpha,\beta(\nu)]-[\alpha(\nu), d^*_A\beta],\label{nohor}
\end{equation}
and note that it holds even if neither $\alpha$ nor $\beta$ is horizontal.

Define a linear map 
\begin{equation*}
T_A: \Lie(\Gk)\to L^2_{k-\frac{5}{2}}(\kb|_{\partial M})
\end{equation*}
given by
\begin{equation}
T_A(g)= d_A\Delta_A g(\nu) + 2(n-1)H_{}\Delta_A g.
\end{equation}
To justify the target space, recall that $\Lie(\Gk)=H^{k+1}_{Dir}(\kb)$.  So, for $g\in\Lie(\Gk)$, we have $d_A\Delta_A g\in H^{k-2}(\kb)$. By Theorem 9.3 in \cite{Palais}, the restriction map from $H^{k-2}(\kb)=L^2_{k-2}\to L^2_{k-\frac{5}{2}}$ is continuous if $k-\frac{5}{2}>0$.  From our assumption that $k>\frac{n}{2}+1$, we have for $n \ge 3$
\begin{equation*}
k-\frac{5}{2} > \frac{n-3}{2}\ge 0.
\end{equation*}
If $n=2$, since $k$ is integral, the condition $k>\frac{2}{2}+1$ forces $k\ge 3$, and thus $k-\frac{5}{2}\ge\frac{1}{2}>0$. Thus, for $n\ge 2$, we see that $T_A$ is well-defined and a continuous operator. Define a set $\mathcal{L}_A\subseteq\Lie(\Gk)$ as 
\begin{equation}
\mathcal{L}_A:=\Span\{\mathcal{R}_A(\alpha,\beta):\alpha,\beta\in H_A\}
\end{equation}
The previous lemma yields
\begin{theorem}\label{bigcor}
The set $\mathcal{L}_A$ is contained in $\mathrm{ker}(T_A)$.  In particular, since 
$T_A$ is not identically $0$, we have that $\overline{\mathcal{L}_A}$ is a proper 
subset of $\Lie(\Gk)$ (where the closure is taken in the $H^{k+1}$ norm).
\end{theorem}
\begin{proof}
For $g\in\mathcal{L}_A$, we have $g\in\mathrm{ker}(T_A)$ by the equation \eqref{curveq} and Lemma \ref{bct}.
\end{proof}
This theorem shows that the image of the curvature form of the Coulomb connection \emph{at one fixed connection} $\nabla_A$ can never be dense in the gauge algebra, unlike the case in \cite{NR}.
%
%
%
%
\subsection{The Smooth $\mathcal{R}_A$ and $\mathrm{ker}(T_A)$}\label{2.4}
While Theorem \ref{bigcor} shows that $\overline{\mathcal{L}_A}$ cannot equal $\Lie(\Gk)$, the closure of the algebra generated by $\mathcal{L}_A$ still may equal $\Lie(\Gk)$.  Indeed, equation \eqref{thecondition1} is not closed under brackets as we will show in Proposition \ref{smooth1}.  In investigating the algebra generated by $\mathcal{L}_A$, we will restrict our attention to $C^\infty$ functions.  Since $C^\infty$ functions are dense in our Sobolev spaces, we do not lose much generality in this restriction, although it will give a denseness result rather than a full Sobolev space result.
So our goal will be to show that $\mathcal{L}_A\cap C^\infty$ algebraically generates $\Lie(\Gk)\cap C^\infty$.  If this is the case, then the closure of the algebra generated by $\mathcal{L}_A$ will be all of $\Lie(\Gk)$.  

We will show that $\mathcal{L}_A\cap C^\infty$ does algebraically generate $\Lie(\Gk)\cap C^\infty$, and thus the closure of the algebra generated by $\mathcal{L}_A$ will be all of $\Lie(\Gk)$, in the certain case where $P$ is the trivial bundle 
$\bar{O}\times K\to \bar{O}$ for a bounded open set $O\subseteq\mathbb{R}^n$ with smooth boundary and where the base connection $\nabla_{A_0}$ is the flat connection.  In this case, 
$\Kb$ is isomorphic to $\bar{O}\times K\to \bar{O}$, and $\kb$ is isomorphic to 
$\bar{O}\times\mathfrak{k}\to \bar{O}$.  We can view gauge transformations $g$ as $K$-valued 
functions on $\bar{O}$, gauge algebra elements $\psi$ as $\mathfrak{k}$-valued 
functions, and $\kb$-valued forms as $\mathfrak{k}$-valued forms.

As in Section~\ref{1.1}, we denote the flat connection as $\nabla_0$.  This means we should denote exterior differentiation
by $d_0$, but since $\nabla_0 = d$ (as asserted in Section~\ref{1.1}), we will instead simply use $d$ without a subscript.  Similiarly, 
we denote $d^*_0$ by $d^*$.

Our first step in showing that $\mathcal{L}_0\cap C^\infty$ algebraically generates $\Lie(\Gk)\cap C^\infty$ is to prove a converse to Lemma \ref{bct}; that is, we will show that
\begin{equation*}
\mathrm{ker}(T_0)\cap C^\infty = \mathcal{L}_0\cap C^\infty.
\end{equation*}

To do this, we will consider slightly different sets than $\mathrm{ker}(T_0)$ and $\mathcal{L}_0$.  Consider the operator $\tilde{T}_0:H^{k-1}(\kb)\to L^2_{k-\frac{5}{2}}(\kb|_{\partial M})$
given by
\begin{equation*}
\tilde{T}_0(f)= d f + 2(n-1)H_{} f.
\end{equation*}
Also, consider the set $\tilde{\mathcal{L}}_0$ defined as
\begin{equation}
\tilde{\mathcal{L}}_0:=\Span\{[\alpha\cdot\beta]:\alpha,\beta\in H_0\}.
\end{equation}
If $f=\Delta g$, note that $f\in\mathrm{ker}(\tilde{T}_0)$ if and only if $g\in\mathrm{ker}(T_0)$, and $f\in\tilde{\mathcal{L}}_0$ if and only if $g\in\mathcal{L}_0$ since $\Delta:H^{k+1}_{Dir}\to H^{k-1}$ is an isomorphism. 
Thus, we have $\mathrm{ker}(T_0)\cap C^\infty = \mathcal{L}_0\cap C^\infty$ if and only if $\mathrm{ker}(\tilde{T}_0)\cap C^\infty = \tilde{\mathcal{L}}_0\cap C^\infty$.  We will prove the latter.

First we look at neighborhoods of the boundary of $O$ and show that all the smooth $\Psi$ that 
satisfy the boundary condition of Lemma \ref{bct} are in $\tilde{\mathcal{L}_0}\cap C^\infty$.
%
%
\begin{lemma}\label{lbo}
Let $O\subset\mathbb{R}^n$ be open and bounded with a smooth boundary and suppose $P=\overline{O}\times K\to \overline{O}$. Let $U\subseteq\overline{O}$ be open in the 
subset topology.  Suppose that $U$ includes a part of the boundary $\partial O$, admits the Type B coordinates  
$\{y_1,\ldots,y_n\}$, and
is the cube $(0,\delta)^{n-1}\times[0,\delta)$ under these coordinates. Let $\Psi\in \mathrm{ker}(\tilde{T_0})\cap C^\infty_c(U;\k)$.
Then $\Psi\in\tilde{\mathcal{L}_0}\cap C^\infty_c(U;\k)$. 
\end{lemma}
In what follows, we shorten ``Dirichlet boundary conditions'' to \emph{DBC}. 
Also, viewing $U$ as the cube $(0,\delta)^{n-1}\times[0,\delta)$, a function $f\in C^\infty_c(U)$ has its support contained in $(\epsilon,
\delta-\epsilon)^{n-1}\times [0,\delta-\epsilon)$ for some 
$\epsilon>0$.  The point is that $f$ need not vanish on the boundary $\{y_n=0\}$.  Lastly, the set $C^\infty_c(U;\k))$ above is 
the set of all $\k$-valued smooth functions on $U$ with compact support.
\begin{proof}[of Lemma \ref{lbo}]
Let $\{v_i\}$ be basis of $\k$.  Then we can write 
$\Psi = \sum \psi_i\cdot v_i$.  Since the basis elements are independent and $\Psi\in \mathrm{ker}(\tilde{T_0})$, 
we have that $d\psi_i(\nu)=-2(n-1)H_{}\psi_i$.  Since 
$\mathfrak{k}$ is semisimple, each basis element $v_i$ can be written as a sum 
of commutators $v_i=\sum_{j=1}^{\alpha(i)}[f_j^i,g_j^i]$.  Hence, we can write 
$\Psi$ as
\begin{equation*}
\Psi = \sum_i \sum_{j=1}^{\alpha(i)}\psi_i [f_j^i,g_j^i].
\end{equation*}
So without loss of generality we can assume $\Psi = \psi\cdot[A,B]$, where 
$A,B$ are fixed elements of $\k$ and $\psi\in C^\infty_c(U)$ and
\begin{equation}\label{psieq} 
d\psi(\nu)=-2(n-1)H_{}\psi.
\end{equation}

Coordinatize $U$ using Type B coordinates $\{y_1,\ldots,y_n\}$ under which 
the domain is the cube $(0,\delta)^{n-1}\times[0,\delta)$. Again let 
$b:=\sqrt{\det(h_{ij})}$, where $h_{ij}$ is the metric tensor of our chart.

Choose $\gamma_1,\gamma_2,\gamma_3,\gamma_4,\gamma_5,\gamma_6\in\mathbb{R}$ so that $supp(\psi)\subset(\gamma_4,\gamma_5)^{n-1}
\times[0,\gamma_5)$ and 
$0<\gamma_1<\gamma_2<\gamma_3<\gamma_4<\gamma_5<\gamma_6<\delta$.  We define a function $\phi$ by setting 
\begin{eqnarray*}
\phi(y_1,\ldots,y_n) &:=& -I(y_1,\ldots,y_{n-2}, y_n)\eta(y_{n-1}) + \\
& & h^{nn}(y_1,\ldots,y_n)b(y_1,\ldots,y_n)^2\psi(y_1,\ldots,y_n),
\end{eqnarray*}
where $\eta$ is a bump function with $\eta\in C_c^{\infty}(\gamma_2,\gamma_3)$ and $\int_{\gamma_2}^{\gamma_3} \eta(s) ds =1$, and $I$ is defined by  
\begin{eqnarray*}
I(y_1,\ldots,y_{n-2}, y_n)&=&\int_{\gamma_4}^{\gamma_5} h^{nn}(y_1,\ldots,y_{n-2},s, y_n)b(y_1,\ldots,y_{n-2},s, y_n)^2\\
& & \psi(y_1,\ldots,y_{n-2},s, y_n)ds.
\end{eqnarray*}  
Then $\phi$ is smooth with compact support.  However $h$ also has an additional property.  Using Proposition \ref{tauy}, Equation \eqref{psieq}, and noting $h^{nn}={1}/{h_{nn}}$, we have
\begin{eqnarray*}
\frac{\partial (h^{nn}b^2\cdot\psi)}{\partial y_n} &=& 
2\sqrt{h^{nn}}\frac{\partial\sqrt{h^{nn}}}{\partial y_n}b^2\psi +
h^{nn} 2b\cdot\frac{\partial b}{\partial y_n}\cdot\psi +
h^{nn}b^2\cdot\frac{\partial\psi}{\partial y_n}\\
&=& h^{nn}b^2\left[2\left(\frac{1}{\sqrt{h^{nn}}}\frac{\partial\sqrt{h^{nn}}}
{\partial y_n} + \frac{1}{b}\cdot\frac{\partial b}{\partial y_n}\cdot
\right)\psi + \frac{\partial\psi}{\partial y_n}\right] \\
&=& 0 \mbox{  on $\partial O\cap U$}.
\end{eqnarray*}
Hence, differentiating under the integral yields
\begin{equation}\label{bco4}
\frac{\partial \phi}{\partial y_n} = 0\mbox{  on $\partial O\cap U$}.
\end{equation}

Define $F:U\to \mathbb{R}$ as
\begin{equation*} 
F(y_1,\ldots,y_n)=\int_0^{y_{n-1}} \phi(y_1,\ldots,y_{n-2},s,y_n)ds.
\end{equation*}  
Then $F$ is smooth and $supp(F)\subset(\gamma_4,\gamma_5)\times(\gamma_1,\gamma_5)\times[0,\gamma_5)$ 
by our construction of $\phi$.  In particular, the term $-I(y_1,\ldots,y_{n-2}, y_n)\eta(y_{n-1})$ was included in the definition of $\phi$ to make $F$ have compact support in the $y_{n-1}$ variable.\footnote{I thank my advisor, Prof. Leonard Gross, for his ideas in making the integral function $F$ have compact support.}  Also, note that 
$F_{n-1} = \phi$ (where the subscript $n-1$ denotes we are taking the partial derivative of $F$ with respect to $y_{n-1}$).  Also, by (\ref{bco4}), differentiating under the integral sign 
yields
\begin{equation}\label{bco5}
F_n = 0\mbox{  on $\partial O\cap U$}.
\end{equation}

We now construct another function $G:[0,\delta]^n\to\mathbb{R}$ by setting 
\begin{equation*}
G(y_1,\ldots,y_n)=\prod_{i=1}^n v_i(y_i), 
\end{equation*}
where each $v_i:[0,\delta]\to\mathbb{R}$ is constructed as follows: for $i \le n-2$,
$v_i\in C^\infty_c(0,\delta)$, $v_i|_{[\gamma_4,\gamma_5]}\equiv 1$, and $supp(v_i)\subset (\gamma_3,\gamma_6)$;
$v_{n-1}\in C^\infty_c(0,\delta)$, $v_{n-1}|_{[\gamma_4,\gamma_5]}(x)=x$, and $supp(v_{n-1})\subset (\gamma_3,\gamma_6)$; 
$v_n\in C^\infty_c([0,\delta))$, $v_n|_{[0,\gamma_5]}\equiv 1$, and $supp(v_n)\subset 
[0,\gamma_6)$. Then
\begin{equation*} 
G_{n-1}|_{(\gamma_4,\gamma_5)^{n-1}\times[0,\gamma_5]}\equiv 1,
\end{equation*} 
and has compact support in $U$.  One can verify that the support of the product $F_{n-1}\cdot G_{n-1}$ lies in $(\gamma_4,\gamma_5)^{n-1}
\times[0,\gamma_5)$, just like the support of $\psi$.  
Furthermore, we have on $(\gamma_4,\gamma_5)^{n-1}\times[0,\gamma_5)$
\begin{equation*}
F_{n-1}\cdot G_{n-1}=h^{nn}b^2\psi,
\end{equation*}
and thus the equation holds everywhere.
Now define 2-forms $\omega_1, \omega_2$ as
\begin{eqnarray*} 
\omega_1 &=& F\cdot A(*^{-1}(dy_1\wedge dy_2\wedge\ldots\wedge dy_{n-2})) \\ 
\omega_2 &=& G\cdot B (*^{-1}(dy_1\wedge dy_2\wedge\ldots\wedge dy_{n-2}))
\end{eqnarray*}
for $n\ge 3$, and 
\begin{eqnarray*} 
\omega_1 &=& F\cdot A(*^{-1}(b^{-1}) \\ 
\omega_2 &=& G\cdot B (*^{-1}(b^{-1}))
\end{eqnarray*}
for $n=2$.
Let $\alpha:=d^*\omega_1$ and 
$\beta:=d^*\omega_2$.  Since $(d^*)^2=0$, we have $d^*\alpha = d^*\beta =0$.  
One can check that 
\begin{equation*}
*(dy_1\wedge\ldots\wedge\widehat{dy_j}\wedge\ldots dy_n) = (-1)^{n+j}\frac{1}{b}\sum_{i=1}^nh_{ij}dx_i.
\end{equation*} 
We have for $n\ge 2$
\begin{eqnarray*}
\alpha &=& d^*(\omega_1) = (-1)^{2n+n+1}*d*(\omega_1) \\
&=& (-1)^{n+1}*d(F\cdot Ady_1\wedge dy_2\wedge\ldots\wedge dy_{n-2})\\
&=& (-1)^{n+1}*((-1)^{n-2}F_{n-1} A dy_1\wedge dy_2\wedge\ldots\wedge dy_{n-2}\wedge dy_{n-1} +\\
& & (-1)^{n-2} F_n A dy_1\wedge dy_2\wedge\ldots\wedge dy_{n-2}\wedge dy_{n}) \\
&=& -\frac{1}{b}\cdot\sum_{i=1}^n (F_{n-1}h_{in} - F_n h_{i(n-1)})dy_iA.
\end{eqnarray*}
Note that in our Type B coordinates we have $h_{in}=0$ 
for $i<n$ everywhere. 
So since $F_n = 0$ on the boundary by (\ref{bco5}), $\alpha$ satisfies DBC 
Similarly,
\begin{equation*}
\beta = -\frac{1}{b}\cdot\sum_{j=1}^n (G_{n-1}h_{jn} - G_n h_{j(n-1)})dy_jB
\end{equation*}
and $\beta$ also satisfies DBC. Indeed, as above $h_{jn}=0$ for $j<n$ everywhere.  Also, since $v_n(y_n)$ is constant on $[0,\gamma_5]$, we have $\left.G_n\right|_{(0,\delta)^{n-1}\times[0,\gamma_5)}=0$, and thus in particular $\left.G_n\right|_{\partial O}=0$.

To calculate $[\alpha\cdot\beta]$, we first note that by the definition of matrix inverses, we have
\begin{equation}\label{inverses}
\sum_{j=1}^n h^{kj}h_{ji} = \delta_{ik}.
\end{equation}
Using the above, we have
\begin{eqnarray*}
[\alpha\cdot\beta] &=& \frac{1}{b^2}\left(\sum_{i,j=1}^n(F_{n-1}h_{in} - F_nh_{i(n-1)})\cdot(G_{n-1}h_{jn}- G_nh_{j(n-1)})h^{ij}\right)[A,B] \\
&=& \frac{1}{b^2}\left(\sum_{i,j=1}^n F_{n-1}G_{n-1}h_{in}h_{jn}h^{ij}-F_{n-1}G_n h_{in}h_{j(n-1)}h^{ij} -\right.\\
& & F_n G_{n-1}h_{i(n-1)}h_{jn}h^{ij}+ F_n G_n h_{i(n-1)}h_{j(n-1)}h^{ij}\Bigg)[A,B] \\
&=& \frac{1}{b^2}\left(\sum_{i=1}^n F_{n-1}G_{n-1}h_{in}\delta_{in}-F_{n-1}G_n h_{in}\delta_{i(n-1)} -\right.\\
& & F_n G_{n-1}h_{i(n-1)}\delta_{in}+ F_n G_n h_{i(n-1)}\delta_{i(n-1)}\bigg)[A,B] \\
&=& \frac{1}{b^2}\left( F_{n-1}G_{n-1}h_{nn}-F_{n-1}G_n h_{(n-1)n} -\right.\\
& & F_n G_{n-1}h_{n(n-1)}+ F_n G_n h_{(n-1)(n-1)}\bigg)[A,B] \\
&=& \frac{1}{b^2}\left( F_{n-1}G_{n-1}h_{nn}+ F_n G_n h_{(n-1)(n-1)}\right)[A,B],
\end{eqnarray*}
where the last line is justified by the fact that $h_{(n-1)n}=0$ 
everywhere in Type B coordinates.  As noted previously, we have $\left.G_n\right|_{(0,\delta)^{n-1}\times
[0,\gamma_5)} = 0$.  Since $supp(F)\subset(\gamma_4,\gamma_5)\times(\gamma_1,\gamma_5)\times[0,\gamma_5)$, we have 
$F_n|_{(0,\delta)^{n-1}\times[\gamma_5,1]} = 0$.  Hence, $F_nG_n\equiv 0$.  So, 
continuing the above, we have
\begin{eqnarray*}
[\alpha\cdot\beta] &=& \frac{1}{b^2}(F_{n-1}G_{n-1}h_{nn})[A,B] \\
&=& \frac{1}{b^2}(h^{nn}b^2\psi h_{nn})[A,B] \\
&=& \psi[A,B],
\end{eqnarray*}
where on the last line we used the fact that $h_{kn}=h^{kn}=0$ for Type B coordinates and thus $h_{nn}h^{nn}=1$ by \eqref{inverses}.  In sum, $\alpha,\beta\in H_0$, and $[\alpha\cdot\beta]=\psi[A,B]$, proving $\Psi=\psi[A,B]\in \tilde{\mathcal{L}}_0\cap C^\infty$, as desired.
\end{proof}
%
%
Next we check that the previous result holds for functions $\Psi$ with compact support.
\begin{proposition}\label{noboundary2}
Let $O\subset\mathbb{R}^n$ be a bounded open set with a smooth boundary and suppose $P=\overline{O}\times K\to \overline{O}$.  Then $C^\infty_c(O;\k)=\Span\{[\alpha\cdot\beta]: \alpha,\beta\in C^\infty_c(\Lambda^1(O;\k)), d^*\alpha =d^*\beta = 0\}\subset\tilde{\mathcal{L}}_0\cap C^\infty$.
\end{proposition}
\begin{proof}
Given $f\in C^\infty_c(O;\k)$, one can cover the support of $f$ with finitely many cubes, and reconstruct $f$ as a product $[\alpha\cdot\beta]$ on each cube in a fashion similar to the process of Lemma \ref{lbo}.  The construction here is simpler since boundary conditions do not come into play.  In particular, one can use the standard coordinates of $\mathbb{R}^n$ whose metric tensor $\{g_{ij}\}$ is of course the identity matrix, greatly simplifying the work. Details for the $n=3$ case can be found in \cite{GrycThesis}.
\end{proof}
We now combine Lemma \ref{lbo} and Propostion \ref{noboundary2} to get our desired global result.
\begin{lemma}\label{curve_im}
Let $O\subset\mathbb{R}^n$ be a bounded open set with a smooth boundary and suppose $P=\overline{O}\times K\to \overline{O}$. Then $\mathrm{ker}(\tilde{T_0})\cap C^\infty=\tilde{\mathcal{L}}_0\cap C^\infty$, and thus $\mathrm{ker}(T_0)\cap C^\infty=\mathcal{L}_0\cap C^\infty$.
\end{lemma}
\begin{proof}
The backward direction has already been shown in Lemma \ref{bct}.  For the 
forward direction, suppose $f\in\mathrm{ker}(\tilde{T_0})\cap C^\infty$, and thus satisfies $df(\nu) = -2(n-1)H_{} f \mbox {  on 
$\partial O$}$.  There exists a finite cover $\{U_k\}_{k=0}^m$ of $\overline{O}$ 
that satisfies the following: $\overline{U}_0\subseteq O$, $\{U_k\}_{k=1}^m$ covers the boundary $\partial O$ and each 
$U_k$ for $k\ge 1$ is a cube in Type B coordinates, and there is a 
partition of unity $\{\lambda_k\}_{k=0}^m$ subordinate to $\{U_k\}_{k=0}^m$ so 
that $d\lambda_k(\nu)=0$ on the boundary.  A construction of such a partition of unity 
can be found in \cite{GrycThesis} for $n=3$ which can easily be generalized for $n\ge 2$. 

With such a partition of unity, we have 
$d(\lambda_k f)(\nu) = -2(n-1)H_{}\lambda_k f \mbox {  on $\partial O$}$.  So, 
by Lemma \ref{lbo} and Proposition \ref{noboundary2} there exists 
$\{\alpha_i^k\}, \{\beta_i^k\}$ such that each 
$\alpha_i^k,\beta_i^k\in C^\infty_c(\Lambda^1(U_k\otimes\k))$, $d^*\alpha_i^k=
d^*\beta_i^k=0$, $\alpha_i^k,\beta_i^k$ satisfy DBC, and
$\lambda_k\cdot f = \sum_{i=1}^{n}[\alpha_i^k\cdot\beta_i^k]$ on $U_k$.  
Extending the $\alpha_i^k$'s and $\beta_i^k$'s by zero, we have 
 $\lambda_k\cdot f = \sum_{i=1}^{n}[\alpha_i^k\cdot\beta_i^k]\in \tilde{\mathcal{L}}_0\cap C^\infty$.  
Since $\tilde{\mathcal{L}}_0$ is a span, $f=\sum_{k=1}^m (\lambda_k\cdot f)\in\tilde{\mathcal{L}}_0$
also, as desired.  So, $\mathrm{ker}(\tilde{T_0})\cap C^\infty=\tilde{\mathcal{L}}_0\cap C^\infty$, and thus $\mathrm{ker}(T_0)\cap C^\infty=\mathcal{L}_0\cap C^\infty$ by the note in the beginning of this subsection.
\end{proof}
%
%
%
%
\subsection{The Generation of the Smooth Gauge Algebra}\label{2.5}
In this section we will use brackets of the image of the curvature form to get 
every smooth function in $\Lie(\Gk)$ for the special case $P=\overline{O}\times K\to \overline{O}$.  The main tool will be 
Lemma \ref{curve_im}.  The first thing we must do is see how the equation
\begin{equation}\label{mainderivative2}
d(\Delta g)(\nu)=-2(n-1)H_{}\Delta g.
\end{equation}
changes when 
we introduce brackets.  More specifically, recall that if $g\in\mathcal{L}_0$, then 
Lemma \ref{curve_im} says $g$ satisfies \eqref{mainderivative2} above.
We want to know how \eqref{mainderivative2} changes if $g$ is 
replaced by $[g_1,g_2]$, for $g_i\in \mathcal{L}_0$.  Indeed, we have
\begin{proposition}\label{smooth1}
Suppose $g_1,g_2\in\mathcal{L}_0$.  Then we have
\begin{equation}
d(\Delta([g_1,g_2]))(\nu)=-2(n-1)H_{}\Delta[g_1,g_2] + 3[\Delta g_1, d g_2(\nu)] 
+3[d g_1(\nu), \Delta g_2].
\end{equation}
\end{proposition}
We state the above proposition for all elements of $\mathcal{L}_0$, not just the smooth elements because the proposition holds in the general case.  However, the use of the proposition in this paper will be just for the smooth case.
\begin{proof}
First note that
\begin{eqnarray*}
\Delta([g_1,g_2]) = [\Delta g_1,g_2]+[g_1,\Delta g_2]-2[d g_1\cdot d g_2].
\end{eqnarray*}
So we have
\begin{eqnarray}
d(\Delta([g_1,g_2]))(\nu) &=& d([\Delta g_1,g_2]+[g_1,\Delta g_2]-2[d g_1\cdot d g_2])(\nu)\nonumber\\
&=& [d(\Delta g_1)(\nu),g_2] + [\Delta g_1,dg_2(\nu)] + [dg_1(\nu),\Delta g_2] \nonumber\\
& & + [g_1,d(\Delta g_2)(\nu)] - 2d([d g_1\cdot d g_2])(\nu).\label{brack1}
\end{eqnarray}
By Lemma \ref{curve_im}, we have
\begin{equation}\label{brack2}
d(\Delta g_i)(\nu) = -2(n-1)H_{}\Delta g_i .
\end{equation}
By equation \eqref{nohor} which follows the proof of Lemma \ref{bct}, we see that if $\alpha,\beta\in H^k_{Dir}(\kb)$ but  
are not necessarily in $H_A$, we have
\begin{equation}
d_A([\alpha\cdot\beta])(\nu) = -2(n-1)H_{}[\alpha\cdot\beta] - [d_A^*\alpha,\beta(\nu)] - [\alpha(\nu), d_A^*\beta].
\end{equation}
Setting $\alpha=d g_1$ and $\beta=d g_2$ above yields
\begin{equation}\label{brack3}
-2d([d g_1\cdot d g_2])(\nu) = -2(-2(n-1)H_{}[d g_1\cdot d g_2] - [\Delta g_1, d g_2(\nu)] - [d g_1(\nu),\Delta g_2]).
\end{equation}
Inserting \eqref{brack2} and \eqref{brack3} into \eqref{brack1}, we have
\begin{eqnarray*}
d(\Delta([g_1,g_2]))(\nu) &=& -2(n-1)H_{}[\Delta g_1,g_2] + [\Delta g_1,dg_2(\nu)] + [dg_1(\nu),\Delta g_2] + \\
& & - 2(n-1)H_{}[g_1,\Delta g_2] + \\ & &-2(-2(n-1)H_{}[d g_1\cdot d g_2] - [\Delta g_1, d g_2(\nu)] - [d g_1(\nu),\Delta g_2])\\
&=& -2(n-1)H_{}([\Delta g_1,g_2]+[g_1,\Delta g_2]-2[d g_1\cdot d g_2]) + \\
& & 3[\Delta g_1, d g_2(\nu)] + 3[d g_1(\nu),\Delta g_2] \\
&=& -2(n-1)H_{}\Delta([g_1,g_2]) + 3[\Delta g_1, d g_2(\nu)] + 3[d g_1(\nu),\Delta g_2],
\end{eqnarray*}
as desired.
\end{proof}
We will now show that the new term in Proposition \ref{smooth1} is actually very 
general.
\begin{proposition}\label{smooth2}
Let $F$ be a smooth $\mathfrak{k}$-valued function on $\partial O$.  Then 
there exists smooth $\k$-valued functions $g_i,h_i\in\mathcal{L}_0$ such that
\begin{equation*}
d(\Delta(\sum_i [g_i,h_i]))(\nu)+2(n-1)H_{}\Delta(\sum_i[g_i,h_i])= F.
\end{equation*}
\end{proposition}
\begin{proof}
Since $\mathfrak{k}$ is semi-simple, there exists $A_i,B_i,C_i\in
\mathfrak{k}$ that  
\begin{equation*}
F =\sum_i f_i[[A_i,B_i],C_i]
\end{equation*}
for some real valued smooth functions $f_i$.  So, without loss of generality, 
assume that $F=f[[A,B],C]$ for some $A,B,C\in\k$.

Take any non-negative, nonzero, real-valued $\phi\in C^\infty_c(O)$.  By the Strong Maximum Principle, we have $G\phi > 0$ for the interior of each connected component of $O$, and thus on the whole interior of $O$ (note that our definition of the Laplacian as $\Delta=d^*d$ means that in local coordinates $\Delta=-\sum_i \frac{\partial}{\partial x_i}$).  Thus, we can apply Lemma 3.4 of \cite{GT} to get 
\begin{equation*}
\frac{\partial (G\phi)}{\partial \nu} > 0. 
\end{equation*}
In particular, $d(G\phi)(\nu)$ never vanishes.  We set $h:= G\phi\cdot C$.  Since $\Delta h = \phi\cdot C$ has 
compact support, $h\in\mathcal{L}_0$ by Lemma \ref{curve_im}.

Let $\{U_k\}_{k=0}^{m}$ be an open cover of $O$ such that $\{U_k\}_{k=1}^{m}$ covers $\partial O$ and $U_k$ are cubes 
in Type A coordinates for $k\ge 1$.  Let $\{\lambda_k\}_{k=1}^m$ be the corresponding partition of unity for the 
cover $\{U_k\cap\partial O\}_{k=1}^m$ of the boundary.  
We set
\begin{equation*}
f_k:=\lambda_k\cdot \frac{f}{3d(G\phi)(\nu)}.
\end{equation*}
In the cube of $U_k$, suppose the $x_n$ interval is $[0,a]$.  Choose a $C^
\infty$ function $\eta:[0,a]\to[0,1]$ such that $\eta|_{[0,a/4]}\equiv 
1$ and $supp(\eta)\subseteq([0,a/2])$. We can extend $f_k$ to a function 
$\tilde{f}$ on $U_k$ by
\begin{equation*}
\tilde{f}(x_1,\ldots,x_n) = f_k(x_1,\ldots,x_{n-1})\eta(x_n)\exp(-2(n-1)H_{}(x_1,\ldots,x_{n-1})x_n).
\end{equation*}
Note that the support of $\tilde{f}$ lies in $U_k$, so $\tilde{f}$ is a function on all of $\overline{O}$.  On $U_k$, we 
have
\begin{eqnarray*}
d\tilde{f}(\nu) &=& \frac{\partial}{\partial x_n}\bigg|_{x_n=0}f_k(x_1,\ldots,x_{n-1})\eta(x_n)\exp(-2(n-1)H_{}(x_1,\ldots,x_{n-1})x_n) \\
&=& -2(n-1)H_{}(x_1,\ldots,x_{n-1})f_k(x_1,\ldots,x_{n-1})\eta(x_n)\cdot \\
& & \exp(-2(n-1)H_{}(x_1,\ldots,x_{n-1})x_n)\\
&=& -2(n-1)H_{} \tilde{f}.
\end{eqnarray*}
By Lemma \ref{curve_im}, the above shows that $G\tilde{f_k}[A,B]\in\mathcal{L}_0$.  Let $g=\sum_{k=1}^m G\tilde{f_k}[A,B]$.
We now verify that $g$ and $h$ were well-chosen.  By Proposition \ref{smooth1} and since $\Delta h|_{\partial O}=\phi C|_{\partial O}\equiv 0$,
\begin{eqnarray*}
d(\Delta([g,h]))(\nu) + 2(n-1)H_{}\Delta[g,h] &=& 3[\Delta g, d h(\nu)] 
+3[d g(\nu), \Delta h] \\
&=& 3[\Delta g, d h(\nu)] \\
&=& 3[\sum_k f_k[A,B], dG\phi(\nu)\cdot C] \\
&=& 3[(\sum_k \lambda_k)\frac{f}{3dG\phi(\nu)}[A,B],dG\phi(\nu)C] \\
&=& f[[A,B],C]]=F,
\end{eqnarray*}
proving the proposition.
\end{proof}
We are now at the point where we can prove the key lemma for our main theorem.  Let $\mathcal{F}$ be the linear space spanned by $\mathcal{L}_0$ and $[\mathcal{L}_0,\mathcal{L}_0]$.
\begin{lemma}\label{mainthm}
Suppose our principal bundle is $P=\overline{O}\times K\to \overline{O}$, where $O\subseteq\mathbb{R}^n$ for $n\ge 2$ is open, bounded and has smooth boundary.  Suppose $g\in\mathrm{Lie}(\Gk)$ and is $C^\infty$.  Then $g\in\mathcal{F}\cap C^\infty$.
\end{lemma}
\begin{proof}
Let $g\in\mathrm{Lie}(\Gk)\cap C^\infty$.
Recall our linear map $T_0:C^\infty(O;\mathfrak{k})\to C^\infty(\partial O;\mathfrak{k})$ defined as
\begin{equation*}
T_0(f) = d(\Delta f)(\nu) + 2(n-1)H_{}\Delta f.
\end{equation*}
Set $u:=T_0(g)$.  By Proposition \ref{smooth2}, there exists a smooth function $f\in\mathcal{F}$ such that $T_0(f)=u$.  Since $T_0$ is 
linear, we have that $T_0(g-f)=0$.  By Lemma \ref{curve_im}, we know that $g-f\in\mathrm{Span(Im}(\mathcal{R}_0))\subseteq\mathcal{F}$.  Hence, $g=f+(g-f)\in\mathcal{F}$, as we desired.
\end{proof}
The preceeding lemma gives us our main result.
\begin{theorem}\label{lastcor}
Consider the trivial principal bundle $\overline{O}\times K\to \overline{O}$, where $O\subseteq\mathbb{R}^n$ for $n\ge 2$ is open, bounded, and has smooth boundary.  Let $k>n/2+1$, and 
suppose $\nabla_{A_0}=\nabla_0$.  The restricted holonomy group $(\Hko)^0(\nabla_0)$ with base point $\nabla_0$ of the Coulomb connection of the associated bundle $\Ck\to\Ck/\Gk$ is dense in the 
connected component of the identity of $\Gk$.
\end{theorem}
Before we prove this theorem, we should mention what we mean by ``holonomy group.''  We define $\Hko(\nabla_0)$ the 
the same way it would be definied in finite dimensions.  That is $g\in\Hko(\nabla_0)$ if and only if $\nabla_0\cdot g$ can be 
connected to $\nabla_0$ by a horizontal path in $\Ck$.  It has been shown that with this definition, $\Hko(\nabla_0)$ is a Banach Lie group, and the restricted holonomy group $(\Hko)^0(\nabla_0)$ is also a Banach Lie group (for the statement of this theorem, see \cite{Vas}).
\begin{proof}
This follows directly from Lemma 7.6 and Proposition 7.7 in \cite{NR}.  Specifically, Lemma 7.6 and the beginning of the proof 
of Proposition 7.7 of \cite{NR} imply that every element of $\mathcal{F}$ is the tangent vector to a curve in $(\Hko)^0(\nabla_0)$.  Then Proposition 7.7 of \cite{NR} tells us that $(\Hko)^0(\nabla_0)$ is dense in the connected component of the identity of $\Gk$ since $\mathcal{F}$ is dense in $\mathrm{Lie}(\Gk)$, completing the proof.
\end{proof}
In sum, we used the image of the curvature $\mathcal{R}$ of the Coulomb connection to tell us about the Lie algebra of 
the holonomy group $\Hko$.  The fact that this image generates the entire holonomy group is a well-known 
theorem in finite dimensions.  A version of this theorem also holds in the infinite-dimensional case, as proved in \cite{Magnot}.
\bibliographystyle{plain}
\bibliography{Will-ref}

\end{document}